\documentclass[12pt, twoside, a4paper]{amsart}
\usepackage{amsfonts, amsmath, amssymb, amsthm, mathtools}
\usepackage{comment}
\usepackage{hyperref}
\usepackage{xcolor}
\hypersetup{colorlinks=true, linkcolor=blue}
\usepackage{hyphenat}
\newtheorem{theorem}{Theorem}
\newtheorem{proposition}{Proposition}

\newtheorem{remark}{Remark}
\newtheorem{lemma}{Lemma}

\newcommand{\R}{\mathbb{R}}

\numberwithin{equation}{section}
\newcommand{\fin}{\hfill\rule{2mm}{2mm}\medskip }
\usepackage[left=1.25in, right=0.75in, top=1.25in, bottom=1 in, includefoot, headheight=14pt]{geometry}

\subjclass[2010]{35R30, 35K10, 35K57, 	93B07 }

\keywords{Parabolic systems, Carleman estimates, Inverse problems, Source estimates.}

\begin{document}
\title[$ L^q$-Carleman estimates and applications to inverse problems]{$ L^q$ Carleman estimates with boundary observations and applications to inverse problems}
\maketitle
{\footnotesize
	
	\centerline{\scshape  Elena-Alexandra Melnig}
	\medskip
	{\footnotesize
		
		\centerline{ Faculty of Mathematics, "Al. I. Cuza" University of Ia\c si, Romania}
		\centerline{Octav Mayer Institute of Mathematics, Romanian Academy, Ia\c si Branch}
		\centerline{\textit{E-mail address}: alexandra.melnig@uaic.ro}
	}}

\begin{abstract}

{We consider  coupled linear parabolic systems and we establish estimates  in $L^q$-norm for the sources  from partial boundary observations of the solutions. The main tool is a family of  Carleman estimates in $L^q$-norm with boundary observations.}

\end{abstract}
\maketitle
\date{}

\section*{Introduction}
In this paper we consider linear parabolic systems  coupled in zero order terms and we  obtain estimates in $L^q$ spaces, $q\geq 2$ for the sources in terms of  measurements of the solutions  on a part of the boundary. This research is based on   previous results in $L^2$ and uses the regularity properties of the heat flow. Such an inverse problem is motivated by the fact that the boundary of a domain is more accessible for measurements on the solution. 

This type of problems was considered by Imanuvilov and Yamamoto in \cite{imayam1998} where they  obtained source estimates in $L^2$ norm with  observations on the solution on a subdomain or on the boundary. Their method  employed the coupling of two  $L^2$-Carleman estimates, a technique no longer applicable  in the $L^q $ framework, as explained in \cite{melnig2020}. Furthermore, they require more detailed information from the measurements, as their estimates involve the solution together with its time and space derivatives, which we overcome in the present paper.   

The problem we address here comes as a continuation of the work  done in \cite{melnig2020} where we have established such  $L^q$ stability estimates for the sources of a linear system in terms of the solution measured on a subdomain. Also, in \cite{melnig2021} we established estimates of this type for the sources in a  semilinear system with observations on the solution in a subdomain. 

An important aspect in the above-mentioned papers as well as in the present study is the positivity of the solutions and of the sources, as the systems we consider are intended to model reaction-diffusion processes. In this context, we need two fundamental tools: strong invariance principles and, respectively, Carleman estimates. 
The first tool relies on the strong maximum principle for parabolic equations and on strong invariance results for weakly coupled parabolic  systems. For such maximum principles or invariance results in the framework of  classical solutions  we refer to  \cite{wein2} for equations and to  \cite{wein} for   systems. For the case of variational solutions we mention \cite{lefter_melnig2024}.

The other main tool in our approach is a family of Carleman estimates in $L^q$ norm, $q\geq 2$, with boundary observations and general weights of exponential type obtained for parabolic systems with general boundary conditions.  The Carleman inequalities in $ L^2$ were first used in the framework of controllability of the heat equation with controls distributed in a subdomain  by O.Yu Imanuvilov, \cite{furima1996}. Since then they have found applications in obtaining observability inequalities in controllability problems, unique continuation properties and inverse problems. 
The Carleman estimates in $L^q$ are derived through a bootstrap argument, starting from the $L^2$-Carleman estimates  proved in \cite{lefter_melnig2024} and using the regularizing effect of the parabolic equation.  This type of argument can be found in V. Barbu \cite{barbu2000}, J.-M. Coron, S.Guerrero, L.Rosier \cite{corgueros2010}, Fernandez-Cara, Zuazua \cite{fernandez-Zuazua2000}, K. Le Balc'h \cite{balch}, in the context of controllability.

\section{Preliminaries and main results}
 
Let $T>0$ and  $\Omega\subset\R^N$ be  an annular  domain, that is a domain which is  diffeomorphic to $B_2(0)\setminus B_1(0)$,  with boundary  $\partial\Omega$ of class $C^2$. Consider $\Gamma_1\subset\partial\Omega$ be the exterior boundary and  $\Gamma_0=\partial\Omega\setminus\overline\Gamma_1$ be the inner boundary, and let $Q=(0,T)\times\Omega$, $\Sigma_1:=(0,T)\times\Gamma_1$, $ \Sigma_0:=(0,T)\times\Gamma_0$.

In the following, whenever we refer to vector-valued functions from a corresponding Sobolev space, we will write $L^q(Q)$ instead of $[L^q(Q)]^n$,   $W^{2,1}_q(Q)\cap L^\infty(Q)$ instead of $[W^{2,1}_q(Q)\cap L^\infty(Q)]^n$. 
\smallskip

We consider linear parabolic systems coupled in zero order terms with the elliptic part of the operator expressed  in divergence form and, for
each component of the system,  with  general homogeneous boundary conditions (Dirichlet, Neumann or Robin) on each connected component of the boundary: 

\begin{equation}\label{sysinitial}\tag{S}
\left\lbrace
\begin{array}{ll}
D_ty_i-\sum\limits_{j,k=1}^N D_j(a_i^{jk} D_ky_i)+\sum\limits_{k=1}^N b_i^{k} D_ky_i+\sum\limits_{l=1}^Nc_i^ly_l=g_i &(0,T)\times\Omega,
\\
\beta_i(x)\frac{\partial y_i }{\partial n_{A}}+\eta_i(x)y_i=0  &(0,T)\times\partial\Omega, \\
\end{array} \quad i=\overline{1,n}
\right.\end{equation} 

where 
\begin{itemize}
\item[(H1)]  $a_i^{jk}\in W^{1,\infty}(\Omega) $, $b^{k}_i,c_i^l  \in  L^\infty(\Omega)$ and $a_i^{jk}$
 satisfy the ellipticity condition: $$\exists \mu>0 \text{ \textit{s.t.} }\sum_{j,k=1}^N a_i^{jk}(x)\xi_j\xi_k\geq\mu|\xi|^2,\quad\forall\xi\in\R^N,\quad(t,x)\in Q,  i=\overline{1,n};$$
\item[(H2)] the sources are positive, $g_i\in L^q(Q)$, $g_i\geq 0, i=\overline{1,n}; $
\item[(H3)] the coupling coefficients are non-positive $c_i^l\leq 0, i\neq l$;
\item[(H4)] $\beta_i, \eta_i\in L^\infty(\partial\Omega), \, \beta_i, \eta_i\geq 0 \text{ and }\beta_i>0  \text{ or } \beta_i\equiv 0 \text{ and }\eta_i\equiv 1$  on each connected component of $\partial\Omega$.
\end{itemize}

The   boundary observation on the solution is  $\zeta=(\zeta_i(y_i))_{i=\overline{1,n}}$: 
\begin{equation}\label{obs}\tag{O}
\zeta_i(y_i)=\gamma_i(x)\frac{\partial y_i}{\partial n_{A}}+\delta_i(x)y_i, \quad (0,T)\times\Gamma_1, \quad \Gamma_1\subset\partial\Omega,\quad i=\overline{1,n},
\end{equation}
with given $\gamma_i,\delta_i\in L^\infty(\Gamma_1)$. 

We impose an independence condition between the boundary conditions and the observations expressed as:

\begin{itemize}
\item[(H5)] $\begin{vmatrix}
\gamma_i&\delta_i\\
\beta_i& \eta_i
\end{vmatrix}\neq 0 \text{ on } \Gamma_1,\quad  i=\overline{1,n}.
$
\end{itemize}
\medskip

For $2\leq q\leq \infty$, $\tilde{c}>0$, $\tilde\delta>0 $ and  $\tilde G$   a compact subset of $ [L^{q'}(Q)]^n$ with $q'=\frac{q}{q-1}$  such that $0\notin\tilde G$, we consider the following classes of sources:

 \begin{equation}\label{Gset}
 \mathcal{G}_{q,\tilde{\delta},\tilde G}=\left.
 \begin{cases}
 &g\in W^{1,1}((0,T);[L^q(\Omega)]^n):\,g\ge 0 \\
 &\text{ and } \exists \tilde{g}\in\tilde G \text{ s.t. } \int_{Q}g\cdot\tilde{g}dxdt\geq\tilde{\delta}\|g\|_{L^q(Q)}
 \end{cases}
 \right\rbrace .
 \end{equation}

 In this context, the main result that gives $L^q$ estimates for the sources assuming that they belong to some class $\mathcal{G}_{q,\tilde{\delta},\tilde G}$ is the following:
  
   \begin{theorem}\label{ThLq}
Consider the system \eqref{sysinitial} with hypotheses $(H1)-(H5)$.  Then for $2\le q<\infty$,    $g\in \mathcal{G}_{q,\tilde{\delta},\tilde G}$ and the corresponding solution $y\in W^{2,1}_q(Q)$,   there exists $C=C(q,\tilde{\delta},\tilde G )>0$ such that 
\begin{equation}\label{gest_q}
	\left\|g\right\|_{L^q(Q)}\leq C_1\|\zeta\|_{L^2(\Sigma_1)}.
	\end{equation} 
	Concerning the $L^\infty$ estimates,  for sources  $g\in \mathcal{G}_{q,\tilde{\delta},\tilde G}$ and corresponding solutions $y\in W^{2,1}_q(Q)$ for all $q<\infty$, there exists $C=C(\tilde{\delta},\tilde G)>0$ such that
	\begin{equation}\label{gest_inf}
	\left\|g\right\|_{L^\infty(Q)}\leq C_1\|\zeta\|_{L^2(\Sigma_1)}.
	\end{equation} 
\end{theorem}

 \section{$L^q-L^2$-Carleman estimates with boundary observations}\label{secLqCar}

The key ingredient to prove the above Theorem is a family of $L^q$ Carleman estimates with observations on the outer boundary. To derive these estimates, it is necessary to use the classical mechanism, but with auxiliary functions satisfying supplementary technical properties. These properties, in turn, dictate the special choice of the shape of the domain.

 \smallskip
 
The choice of auxiliary functions is the following:
\begin{equation}
\psi_0\in C^2(\overline{\Omega}) \text{ s.t. } \psi_0|_{\Gamma_0}=0,\quad\psi_0|_{\Gamma_1}=1,\quad |\nabla\psi_0|>0 \text{ in } \overline{\Omega}, \quad\frac{\partial\psi_0}{\partial \nu}\left|_{\Gamma_0<0}\right.
\end{equation}
 $\psi=\psi_0+K$, $K$ big enough such that 
$$\frac{\sup\psi}{\inf\psi}\leq\frac{8}{7}$$
 and the weight functions
\begin{equation}
\varphi_0(t,x):=\frac{e^{\lambda\psi_0(x)}}{t(T-t)},\quad
\alpha_0(t,x):=\frac{e^{\lambda\psi_0(x)}-e^{1.5\lambda\|\psi_0\|_{C(\overline{\Omega})}}}{t(T-t)}
\end{equation}
\begin{equation}
\varphi(t,x):=\frac{e^{\lambda\psi(x)}}{t(T-t)}\,,\quad
\alpha(t,x):=\frac{e^{\lambda\psi(x)}-e^{1.5\lambda\|\psi\|_{C(\overline{\Omega})}}}{t(T-t)}. 
\end{equation}

 \begin{remark}
The choice of annular domains is imposed by technical constraints necessary to derive the appropriate Carleman estimates. More specifically, these constraints arise from the conditions imposed on the function $\psi_0$, which require the presence of two connected components of the boundary with the same topological structure. Consequently, while alternative domain choices are possible, they must satisfy the fundamental requirement that a function with the necessary properties exists. 

 \end{remark}
 
The Carleman estimate we obtain is the following.
\begin{proposition}\label{lemaCarlemanLq}
 Let  $g\in L^q(Q)$, with $2\leq q<\infty$.  Then there exist   $s_0=s_0(q)$, $\lambda_0=\lambda_0(q)$,  such that for $\lambda>\lambda_0$, $s', s>s_0$, $\frac{s'}{s}>\bar\gamma>1$,  there exists $C=C(q,\bar\gamma)$:
\begin{equation}\begin{aligned}
& \|ye^{s'\alpha}\|_{L^{q}(Q)}+ \|(Dy)e^{s'\alpha}\|_{L^{q}(Q)}+\|(D^2y)e^{s'\alpha}\|_{L^{q}(Q)}+\|(D_ty)e^{s'\alpha}\|_{L^{q}(Q)}\\
&\leq C\left[\|ge^{s{\alpha}}\|_{L^{q}(Q)}+\|\zeta e^{s\alpha}\|_{L^2(\Sigma_1)}\right]
 \end{aligned}
 \end{equation}
  where the constant  $C$ depends on $\lambda$ but independent of $s$.
 \end{proposition}

 The result in Proposition \ref{lemaCarlemanLq} relies on  the regularizing effect of the parabolic flow combined with a  bootstrap argument applied to the linear  parabolic  system and using  the following  $L^2$ Carleman estimate that was obtained in \cite{lefter_melnig2024}, adapting the classical techniques from \cite{furima1996}:
\begin{proposition}\label{lemaCarlemanclassic} For $g\in L^2(Q)$,
there exist constants $\lambda_0=\lambda_0(\Omega),$ $s_0=s_0(\Omega)$ such that, for any $\lambda\geq\lambda_0$, $ s\geq s_0$ and some $C=C(T,\Omega)$, the following inequality holds:
\begin{equation}\label{classicalCarleman}
\begin{aligned}
&\int_{Q}\left[(s\varphi)^{-1}\left(|D_ty|^2+|D^2y|^2\right)+s\lambda^2\varphi|Dy|^2+s^3\lambda^4\varphi^{3}|y|^2\right]e^{2s\alpha}dxdt+
\int_{[0,T]\times\Gamma_0}s^3\lambda^3\varphi^{3}|y|^2e^{2s\alpha}d\sigma\\
&\leq
C\left(\int_{Q}|g|^2e^{2s\alpha}dxdt+
\int_{[0,T]\times\Gamma_1}s^3\lambda^3\varphi^{3}|\zeta|^2e^{2s\alpha}d\sigma\right),\qquad\quad
\end{aligned}
\end{equation}
for $y\in H^1(0,T; L^2(\Omega))\cap L^2(0,T; H^2(\Omega))$ solution of \eqref{sysinitial}.
\end{proposition}

In the following  we denote the operators entering the system by 

 $B=(B_i)_{i=\overline{1,n}}$, $B_iy_i=\beta_i(x)\frac{\partial y_i}{\partial n_{A_i}}+\eta_i(x)y_i$ , $L=(L_i)_{i=\overline{1,n}}$ , $L_iy_i=-\sum\limits_{j,k=1}^N D_j(a_i^{jk} D_ky_i)$, $L^1=(L^1_i)_{i=\overline{1,n}}$, $L^0=(L^0_i)_{i=\overline{1,n}}$,  where the lower-order operators are given by ($w$ is a scalar function, $y$ is vector valued function):
\begin{equation}\label{secpart} L^1_iw=\sum\limits_{{\substack{k=\overline{1,N}}}}b_i^{k} D_kw,\quad L^0_iy=\sum\limits_{l=\overline{1,n}}c^{l}_i y_l, \quad i=\overline{1,n}.\end{equation}

Using these notations the system \eqref{sysinitial} is written    more compactly

\begin{equation}
\begin{cases}
&D_ty+Ly+L^1y+L^0y=g, \,  (0,T)\times\Omega,\\
& By=0,\, (0,T)\times\partial\Omega.
\end{cases}
\end{equation}

For  the proof of $L^q$ estimates we need  Sobolev embedding results for anisotropic Sobolev spaces ( see \cite{lady}, Lemma 3.3):
\begin{lemma}\label{lemmaLady}
	Consider $u\in W^{2,1}_p(Q)$.
	
	Then $u\in Z_1$ where 
	$$
	Z_1= \left\{\begin{array}{lll}
	L^q(Q)&\text{ with }q\le \frac{(n+2)p}{n+2-2p}&\text{ when }p<\frac{N+2}{2}\\
	L^q(Q)&\text{ with }q\in[1,\infty),&\text{ when }p=\frac{N+2}{2}\\
	C^{\alpha,\alpha/2}(Q)& \text{ with } 0<\alpha <2-\frac{N+2}{p},&\text{ when }p>\frac{N+2}{2}
	\end{array}\right.
	$$	
	and there exists $C=C(Q,p,N)$ such that $$
	\|u\|_{Z_1}\le C\|u\|_{W^{2,1}_p(Q)}.
	$$
	Moreover, $Du\in Z_2$ where
	$$
	Z_2= \left\{\begin{array}{lll}
	L^q(Q)&\text{ with }q\le \frac{(N+2)p}{N+2-p}&\text{ when }p<{N+2}\\
	L^q(Q)&\text{ with }q\in[1,\infty),&\text{ when }p={N+2}\\
	C^{\alpha,\alpha/2}(Q)& \text{ with } 0<\alpha<1-\frac{N+2}{p},&\text{ when }p>{N+2}
	\end{array}\right.
	$$	
	and there exists $C=C(p,N)$ such that $$
	\|Du\|_{Z_2}\le C\|u\|_{W^{2,1}_p(Q)}.$$
\end{lemma}

We   introduce the following auxiliary functions which do not depend on space variable:
$$
\overline{\varphi}:=\frac{e^{\lambda(K+1)}}{t(T-t)},\quad \underline{\varphi}:=\frac{e^{\lambda K}}{t(T-t)},\quad \overline{\alpha}:=\frac{e^{\lambda(K+1)}-e^{1.5\lambda(K+1)}}{t(T-t)},\quad \underline{\alpha}:=\frac{e^{\lambda K}-e^{1.5\lambda(K+1)}}{t(T-t)}.
$$

\begin{remark}\label{remarkweightsbound}
Observe that for some $\sigma>\sigma_0$ and for some $\tilde\lambda_0(\sigma_0)>0$ we have for $\lambda>\tilde\lambda_0$ that
$$
-(\sigma-1)e^{1.5\lambda(K+1)}+\sigma e^{\lambda(K+1)}-e^{\lambda K}\le -\sigma \lambda e^{\lambda(K+1)}.
$$
Which gives 
$$
\displaystyle\frac{e^{m\lambda(K+1)}\sigma^ms_1^m\lambda^m}{t^m(T-t)^m}e^{\frac{-(\sigma-1)s_1e^{1.5\lambda(K+1)}+\sigma s_1e^{\lambda(K+1)}-s_1e^{\lambda K}}{t(T-t)}}\leq \frac{e^{m\lambda(K+1)}\sigma^ms_1^m\lambda^m}{t^m(T-t)^m}  e^\frac{-\sigma s_1\lambda e^{\lambda(K+1)}}{{t(T-t)}}$$
$$ \leq\sup_{\mu\in[0,\infty)}\mu^m e^{-\mu}=C(m),$$
which gives
\begin{equation}\label{2.2}
\overline\varphi^ms_2^m\lambda^me^{s_2\overline\alpha}\leq C(m)e^{s_1\underline\alpha}.
\end{equation}
In conclusion, for all $m>0$  and $\sigma_0>1$,  there exist $\tilde\lambda_0=\tilde\lambda_0(\sigma_0)>0$ and $C=C(m)$ such that if $\lambda>\lambda_0$ and $s_1,s_2>0$ with $\frac{s_2}{s_1}=\sigma>\sigma_0$, one has 
\begin{equation}
\varphi^ms_2^m\lambda^me^{s_2\alpha}\leq C(m)e^{s_1\alpha},
\end{equation}
with $\varphi,\alpha$ as above.
\end{remark}

The auxiliary functions are constructed based on $\psi_0$  following the same methodology used in the case of internal observations. Moreover, the associated weight functions maintain analogous estimates, preserving their  behaviour when working to obtain the Carleman estimates. This allows for a consistent framework in handling both cases and ensures that the results obtained for internal observations extend naturally to the current setting(see  \cite{melnig2021}).

\textbf{Proof of Proposition \ref{lemaCarlemanLq}.\,}
For a given $\gamma>1$ and $j\in \mathbb{N}$  we define 
\begin{equation*}
w^j:=ye^{\gamma^js\overline\alpha}=w^{j-1}e^{\gamma^{j-1}s\overline{\alpha}(\gamma-1)}.
\end{equation*} 
Notice that since $\gamma>1$, it follows that  for a  fixed $j$ there exists $\bar\lambda_0(j)>0$ such that 
\begin{equation}\label{comparatii}
e^{\gamma^{j}s\alpha}<e^{\gamma^{j}s\overline\alpha}<e^{s\alpha}
\text{ for all } \lambda\geq \bar\lambda_0(j).
\end{equation}
Moreover, each   $w^j$ verifies the initial boundary value problem
\begin{equation}\label{sysz}
\left\{\begin{aligned}
&D_tw^j+Lw^j=ge^{\gamma^{j}s\overline{\alpha}}
+O[s\gamma^j\overline\varphi^2e^{\gamma^{j-1}s\overline{\alpha}(\gamma-1)}
]w^{j-1},\text{ in }(0,T)\times\Omega\\
&Bw^j=0,\text{ on } (0,T)\times\partial\Omega,\\
&w^j(0,\cdot)=0 \text{ in }\Omega .
\end{aligned}\right.\quad
\end{equation} 

Observe that the boundary conditions satisfied by $w^j$ remain the same as those satisfied by $y^j$ as we have used in the definition of $w^j$ the weight $\overline\alpha$ which is independent of the space variables.

We describe the bootstrap argument: based on the the regularity argument from  Lemma \ref{lemmaLady},  we construct the sequence $\{q_j\}_{j\in\mathbb{N}}$:
\begin{equation}\label{seqqj}
q_0=2,\quad q_j:=
\begin{cases}
\dfrac{(N+2)q_{j-1}}{N+2-q_{j-1}}, \text{  if } q_{j-1}<N+2,\\
\frac{3}{2}q_{j-1}, \text{  if } q_{j-1}\geq N+2.
\end{cases}
\end{equation}
Observe that the sequence $\{q_j\}_{j\in\mathbb{N}}$ is increasing to infinity. Now, since $g\in L^q(Q)$, we can take $m$ such that  $q_{m-1}\leq q<q_m$ and we get,  by standard Sobolev embedding, that: 

\begin{equation}
\label{sobo} W^{1,q_{j-1}}(Q)\subset L^{q_j}(Q),\, j=1,\ldots,m, 
\end{equation}
because the Sobolev exponent $q_j^*:=\frac{(N+1)q_{j-1}}{N+1-q_{j-1}}$ is greater than $q_j$.

Now,  an argument like the one in  Remark \ref{remarkweightsbound} gives that there exist $S_0,\Lambda_0\geq 0 $ and $C=C(j)>0$ such that for $s\geq S_0, \lambda\geq\Lambda_0$ and 
 we have
\begin{equation}\label{obsbound}
s\gamma^j\overline\varphi^2e^{\gamma^{j-1}s\overline{\alpha}(\gamma-1)}\leq C(j) \text{ for }j=\overline{1,m}.
\end{equation}
Since the initial data $w^j(0,\cdot)$ is zero and using the previous estimate \eqref{obsbound}, the parabolic regularity gives that for $\lambda$ big enough ($\lambda>\max\{\lambda_0,\Lambda_0,\max_{j=\overline{1,m}}\{\tilde\lambda_j,\bar\lambda_j\}\}$) we get
\begin{equation}\label{parabolicregz}
\|w^j\|_{W^{2,1}_{q_{j-1}}(Q)}\leq C\left(\|ge^{\gamma^{j}s\overline{\alpha}}
\|_{L^{q_{j-1}}(Q)}+\|w^{j-1 }\|_{L^{q_{j-1}}(Q)}\right).
\end{equation}
Using Lemma \ref{lemmaLady} and the embedding \eqref{sobo} we get that
$w^j \in L^{{q}_{j}}(Q) $ and, from the previous inequality, we obtain the estimate  
\begin{equation}\label{embedding}
\|w^j\|_{L^{q_{j}}(Q)}\leq C\left(\|ge^{\gamma^{j}s\overline{\alpha}}
\|_{L^{q_{j-1}}(Q)}+\|w^{j-1 }\|_{L^{q_{j-1}}(Q)}\right) \text{ for } j=1,\ldots,m, 
\end{equation}
giving
\begin{equation}\label{embedding1}
\|w^m\|_{L^{q_{m}}(Q)}\leq C\left(\sum_{j=1}^{m-1}\|ge^{\gamma^{j}s\overline{\alpha}}
\|_{L^{q_{j}}(Q)}+\|w^{0}\|_{L^{q_{0}}(Q)}\right).
\end{equation}
Regarding the first order terms, Lemma \ref{lemmaLady} implies  that 
$\|Dw^m\|_{L^{q_{m}}(Q)}\leq \|w^m\|_{W^{2,1}_{q_{m-1}}(Q)}$, and  we  obtain estimates for the first-order derivatives 
 \begin{equation}\label{embedding2}
 \|w^m\|_{L^{q_{m}}(Q)}+\|Dw^m\|_{L^{q_{m}}(Q)}\leq C\left(\sum_{j=1}^{m-1}\|ge^{\gamma^{j}s\overline{\alpha}}
 \|_{L^{q_{j}}(Q)}+\|w^{0}\|_{L^{q_{0}}(Q)}\right).
 \end{equation}
  Since we have 
 $$\|ye^{\gamma^{m}s\alpha}\|_{L^q(Q)}\leq\|ye^{\gamma^{m}s\overline\alpha}\|_{L^q(Q)}\leq C\|ye^{\gamma^{m}s\overline\alpha}\|_{L^{q^m}(Q)} $$
 and 
 $$\|Dye^{\gamma^{m}s\alpha}\|_{L^q(Q)}\leq\|Dye^{\gamma^{m}s\overline\alpha}\|_{L^q(Q)}\leq C\|Dye^{\gamma^{m}s\overline\alpha}\|_{L^{q^m}(Q)},$$
 using \eqref{embedding2} we get an estimate,
  \begin{equation}\label{embedding4}
 \|ye^{\gamma^{m}s\alpha}\|_{L^q(Q)}+\|Dye^{\gamma^{m}s\alpha}\|_{L^q(Q)}\leq C\left(\sum_{j=1}^{m-1}\|ge^{\gamma^{j}s\overline{\alpha}}
 \|_{L^{q_{j}}(Q)}+\|w^{0}\|_{L^{q_{0}}(Q)}\right).
 \end{equation}
 
Since ${q_0=2}$ and $w^0=ye^{s\overline{\alpha}}$, using  \eqref{comparatii} and     the $L^2$ Carleman inequality \eqref{classicalCarleman} we can bound $ \|w^{0}\|_{L^{q_{0}}(Q)}$:
\begin{equation}
\|ye^{s\overline{\alpha}}\|_{L^2(Q)}\leq C\|ye^{\frac{1}{\gamma}s\alpha}\|_{L^2(Q)}\leq C\left(\|ge^{\frac{1}{\gamma}s\alpha}\|_{L^2(Q)}+\|s^{\frac{3}{2}}\lambda^2\varphi^{\frac{3}{2}}\zeta e^{\frac{1}{\gamma}s\alpha}\|_{L^2(\Sigma_1)}\right).
\end{equation}
Now, since $q>q_j$ forall $ q_j\in\overline{1,m-1}$, we may use \eqref{comparatii} and  Remark \ref{remarkweightsbound} to get that  there exists $C>0$ such that the right hand-side of \eqref{embedding2} obeys the estimate
\begin{equation}\label{rightbound}\begin{aligned}
&\sum_{j=1}^{m-1}\|ge^{\gamma^{j}s\overline{\alpha}}
 \|_{L^{q_{j}}(Q)}+\|ge^{\frac{1}{\gamma}s\alpha}\|_{L^2(Q)}+\|s^{\frac{3}{2}}\lambda^2\varphi^{\frac{3}{2}}\zeta e^{\frac{1}{\gamma}s\alpha}\|_{L^2(\Sigma_1)}\\
 &\leq C\left(\|ge^{\frac{1}{\gamma}s\alpha}\|_{L^q(Q)}+\|\zeta e^{\frac{1}{\gamma^2}s\alpha}\|_{L^2(\Sigma_1)}\right).\\
\end{aligned}
\end{equation}
To estimate the time derivatives and second-order spatial derivatives,  we  look at the  problem verified by $w^{m+1}$:
\begin{equation}\label{syszm1}
\left\{\begin{aligned}
&D_tw^{m+1}+Lw^{m+1}=ge^{\gamma^{m+1}s\overline{\alpha}}
+O[s\gamma^{m+1}\overline\varphi^2e^{\gamma^{m}s\overline{\alpha}(\gamma-1)}
]w^{m},\text{ in }(0,T)\times\Omega\\
&Bw^{m+1}=0,\text{ on } (0,T)\times\partial\Omega,\\
&w^{m+1}(0,\cdot)=0 \text{ in }\Omega .
\end{aligned}\right.\quad
\end{equation} 

Results of parabolic regularity  together with equation  \eqref{obsbound} yields
\begin{equation}\label{parabregzm1}
\|w^{m+1}\|_{W^{2,1}_{q}(Q)}\leq C\left(\|ge^{\gamma^{m+1}s\overline{\alpha}(\gamma-1)}
\|_{L^{q}(Q)}+\|w^{m}\|_{L^{q}(Q)}\right).
\end{equation}
This implies, once again using \eqref{comparatii}, \eqref{embedding2}, and \eqref{rightbound}, that
\begin{equation}\label{D2ybound}
\begin{aligned}
&\|D^2ye^{\gamma^{m+1}s\alpha}\|_{L^q(Q)}\leq\|D^2ye^{\gamma^{m+1}s\overline\alpha}\|_{L^q(Q)}=\|D^2w^{m+1}\|_{L^{q}(Q)}\\
&\leq C\left(\|ge^{\frac{1}{\gamma}s\alpha}\|_{L^q(Q)}+\|\zeta e^{\frac{1}{\gamma^2}s\alpha}\|_{L^2(\Sigma_1)}\right)
\end{aligned}\end{equation}
and using \eqref{obsbound}, \eqref{comparatii}, \eqref{embedding2}, \eqref{rightbound} that

\begin{equation}\label{Dtybound}
\begin{aligned}
&\|(D_ty)e^{\gamma^{m+1}s\alpha}\|_{L^q(Q)}\leq\|D_tw^{m+1}\|_{L^q(Q)}+C\|w^{m}\|_{L^{q}(Q)}\\
&\leq C\left(\|ge^{\frac{1}{\gamma}s\alpha}\|_{L^q(Q)}+\|\zeta e^{\frac{1}{\gamma^2}s\alpha}\|_{L^2(\Sigma_1)}\right).
\end{aligned}
\end{equation}
Now, using \eqref{embedding2},\eqref{rightbound},\eqref{D2ybound} and \eqref{Dtybound}, for  $\gamma=\bar\gamma^\frac{1}{m+3}$ and  $s$ changed into $\frac{1}{\gamma^2}s$, we conclude that
\begin{equation}\begin{aligned}
& \|ye^{s'\alpha}\|_{L^{q}(Q)}+ \|(Dy)e^{s'\alpha}\|_{L^{q}(Q)}+\|(D^2y)e^{s'\alpha}\|_{L^{q}(Q)}+\|(D_ty)e^{s'\alpha}\|_{L^{q}(Q)}\le\\
&\leq C\left(\|ge^{s{\alpha}}\|_{L^{q}(Q)}+\|\zeta e^{s\alpha}\|_{L^2(\Sigma_1)}\right).
 \end{aligned}\end{equation}
 
 \begin{remark} Observe that if $q>N+1$, then the Morrey embedding theorem gives $L^\infty$ estimates for $y$ and $Dy$: there exist   $s_0,\lambda_0>0$,  s.t. for $\lambda>\lambda_0$, $s', s>s_0$, $\frac{s'}{s}>\bar\gamma>1$,  there exists $C=C(\bar\gamma)$:
 \begin{equation}\label{embedding5}
 \|ye^{s'\alpha}\|_{L^\infty(Q)}+\|Dye^{s'\alpha}\|_{L^\infty(Q)}\leq C\left(\|ge^{s\alpha}
 \|_{L^{q}(Q)}+\|\zeta e^{s\alpha}\|_{L^{2}(\Sigma_1)}\right).
 \end{equation}
 \end{remark}
\fin

\section{Source stability for linear systems. Proof of Theorem \ref{ThLq}}

 In the following, the strong solutions $y\in W^{2,1}_q(Q)$ that we work with are also variational solutions of the system \eqref{sysinitial}. 
 
 For this purpose, we denote by $\Gamma_D^i=\{ x\in\partial\Omega | \beta_i(x)=0\}$ the Dirichlet boundary corresponding to $y_i$. Let us  consider the Hilbert spaces
 $V_i=\{v\in H^1(\Omega):v=0\text{ on }\Gamma_D^i\}$ and let $V=V_1\times\cdots\times V_n$.
\smallskip
 
 Then, for some initial data  $y_0\in L^2(\Omega)$, $y\in L^2(0,T;V)\cap H^1(0,T; V')\subset C([0,T];H)$ is a weak solution of system \eqref{sysinitial} if:
\begin{equation}
\label{varsol1}
\begin{aligned}
&\sum\limits_{i=1}^n(\langle y_i(t),v_i\rangle_{L^2(\Omega)}-\langle y_{0,i},v_i\rangle_{L^2(\Omega)})+\sum\limits_{i=1}^n\int_0^t \langle \mathbb{A}_i\nabla y_i(\tau,\cdot),\nabla v_i\rangle_{L^2(\Omega)}d\tau\\
&+\sum\limits_{i=1}^n\int_0^t \langle b_i\cdot\nabla y_i(\tau,\cdot), v_i\rangle_{L^2(\Omega)}d\tau+\sum\limits_{l=1}^n \int_0^t \langle c_i^ly_l(\tau,\cdot),v_i\rangle_{L^2(\Omega)}d\tau\\
 &+\sum\limits_{i=1}^n\int_0^t\langle \eta_i y_i(\tau,\cdot), v_i \rangle_{L^2(\partial\Omega\setminus\Gamma_{D}^i)} d\tau
 =\sum\limits_{i=1}^n\int_0^t\langle g_i(\tau,\cdot), v_i\rangle_{L^2(\Omega)}d\tau,\quad  \forall t\in (0,T),\quad \forall v_i\in V_i,
\end{aligned}
\end{equation}
where  $\mathbb{A}_i$ are the matrix coefficients matrix $(a_i^{jk})_{j,k\in\overline{1,N}}$ and $b_i=(b_i^k)^\top_{k\in\overline{1,N}}$.

 \textbf{Proof of Theorem \ref{ThLq}}. 
We have to prove that there exists $C=C(q,\tilde{\delta},\tilde G)$ such that for $g\in\mathcal{G}_{q,\tilde{\delta},\tilde G}$, 
\begin{equation}
\begin{aligned}
&\left\|g\right\|_{L^q(Q)}\leq C\|\zeta(y)\|_{L^2(\Sigma_1)}.
\end{aligned}
\end{equation}
We argue by contradiction. Then there  exists a sequence of sources denoted as $(g^m)_m\subset \mathcal{G}_{q,k}$, and  the corresponding solutions $(y^m)_m\subset W^{2,1}_q(Q)$, such that:

\begin{equation}\label{contrad1}
\begin{aligned}
&\left\|g^m\right\|_{L^q(Q)}>m\|\zeta(y^m)\|_{L^2(\Sigma_1)}.
\end{aligned}
\end{equation}
Without loss of generality we can presume that $\left\|g^m\right\|_{L^q(Q)}=1$. For the case when $2\leq q<\infty$, up to a subsequence denoted again $(g^m)_m$, $g^m\rightharpoonup g$ weakly $L^q(Q)$  for some $g\in L^q(Q)$. For the case when $q=\infty$ we have that up to a subsequence  $g^m\rightharpoonup g$ weak-* in $L^\infty(Q)$ for some $g\in L^\infty(Q)$.

This means that the above equation \eqref{contrad1} gives for the sequence of observations that
\begin{equation}\label{convomega}
\zeta(y^m)\rightarrow 0 \text{ in } L^2(\Sigma) \text{ as } m\rightarrow\infty. 
\end{equation}
\smallskip

Observe that in both cases  the weak limit $g$ of $(g^m)_m$ is not zero. Indeed,
for the case when $2\leq q<\infty$,  since $g^m\in\mathcal{G}_{q,\tilde{\delta},\tilde G}$, there exists a corresponding $\tilde{g}^m\in\tilde{G}$ such that
 $$\langle g^m,\tilde{g}^m\rangle_{L^q,L^{q'}}=\int_{Q}g^m\tilde{g}^m\geq\tilde{\delta}\|g^m\|_{L^q(Q)}=\tilde{\delta}.$$
By extracting a further subsequence if necessary, we may assume that \(\tilde{g}_m \to \tilde{g} \in \tilde{G}\), with \(\tilde{g} \neq 0\), strongly in \(L^{q'}(Q)\). Given the weak convergence of \((g^m)_m\) in \(L^q(Q)\) and the strong convergence of \((\tilde{g}^m)_m\) in \(L^{q'}(Q)\), it follows that  

\[
\int_Q g \tilde{g} \geq \tilde{\delta} > 0,
\]  

meaning  that \(g \not\equiv 0\).

Also for the case when $q=\infty$, the weak-* limit $g$ of $(g^m)_m$ is not zero. Since $g^m\in\mathcal{G}_{q,\tilde{\delta},\tilde G}$ and $\tilde{g}^m\in\tilde{G}\subset L^1(Q)$,
 $$\langle g^m,\tilde{g}^m\rangle_{L^\infty,L^{1}}\geq\tilde{\delta}\|g^m\|_{L^\infty(Q)}=\tilde{\delta}.$$
 Because $(g^m)_m$  converges weak-* in $L^\infty(Q)$  we have that $\int_{Q}g\tilde{g}^m\ge\tilde\delta>0$ and thus $g\not\equiv0$.
\medskip

Now, we pass to the limit in the  weak formulation of the problem  \eqref{varsol1} for some initial data  $y_0\in L^2(\Omega)$. 

Consider $(y^m)_m$ a   sequence of solutions for the system \eqref{sysinitial} with corresponding sources $(g^m)_m\subset L^q(Q)$.
For $2\leq q<\infty$, $(g^m)_m$ is bounded  in $L^q(Q)$, $(\zeta^m)_m=(\zeta(y^m))_m$ is bounded in $L^2(\Sigma_1)$ and using the  Carleman estimate for the   solution $y^m$  corresponding to the sources $g^m$, we have that
\begin{equation}
\begin{aligned}
& \|y^me^{s'\alpha}\|_{L^{q}(Q)}+ \|(Dy^m)e^{s'\alpha}\|_{L^{q}(Q)}+\|(D^2y^m)e^{s'\alpha}\|_{L^{q}(Q)}+\|(D_ty^m)e^{s'\alpha}\|_{L^{q}(Q)}\\
&\leq C\left[\|g^me^{s{\alpha}}\|_{L^{q}(Q)}+\|\zeta^m e^{s\alpha}\|_{L^2(\Sigma_1)}\right],
\end{aligned}
 \end{equation}
which gives that $(y^m)_m$ is bounded in ${W^{2,1}_q( Q^\epsilon)}$, for a cylinder of form $ Q^\epsilon=(\epsilon, T-\epsilon)\times\Omega$ with $\epsilon>0$ fixed $0<\epsilon<\frac{T}{2}$, arbitrarily small.

Thus $(y^m)_m$ is bounded in $L^q(\epsilon,T-\epsilon;W^{2,q}(\Omega))$ and $D_t y^m$ is  bounded in $L^q(\epsilon,T-\epsilon;L^{q}(\Omega))$, for $2\leq q<\infty$. We have  \( W^{2,q}(\Omega) \subset W^{1,q}(\Omega) \subset L^q(\Omega) \) with compact embeddings and, in this context,  the Aubin-Lions lemma applied to the sequence \( (y^m)_m \) ensures the existence of a function \( y \in L^q_{\text{loc}}(0,T;W^{1,q}(\Omega)) \) and a subsequence, still denoted \( (y^m)_m \), such that  

\[
y^m \to y \quad \text{strongly in} \quad L^q(\epsilon,T-\epsilon;W^{1,q}(\Omega)) \quad \text{as} \quad m \to \infty, \quad \forall \epsilon > 0.
\]

The strong convergence of \( (y^m)_m \) in \( L^q(\epsilon,T-\epsilon;W^{1,q}(\Omega)) \), together with the weak convergence of \( (g^m)_m \) in \( L^q(Q) \) for \( 2 \leq q < \infty \), allows us to pass to the limit in the variational formulation of problem \eqref{varsol1}. This establishes that \( y \) is the solution to \eqref{sysinitial} corresponding to \( g \). Moreover, in this setting, for all \( \varepsilon > 0 \), we also obtain weak convergence of \( (y^m)_m \) to \( y \) in \( L^q(\epsilon,T-\epsilon;W^{2,q}(\Omega)) \).

When $q=+\infty$ we have, up to a subsequence,  weak-* convergence of $(g^m)_m$ to $g$ in $L^\infty(Q)$  and weak $L^p(Q)$ for all $2\le p<+\infty$. By the above argument we find in fact strong convergence, up to a subsequence, of $(y^m)_m$ to  $y$ in $ L^p(\epsilon,T-\epsilon;W^{1,p}(\Omega))$, as well as weak convergence in $L^p(\epsilon,T-\epsilon;W^{2,p}(\Omega))$, for $2\leq p<\infty$, and so $y$ is solution corresponding to the source $g$.  

Considering that for $2\le p<\infty$ the observation $\zeta$ is linear continuous operator in $L(W^{2,p}(\Omega),L^p(\Gamma_1))$ 
we have that 
$$\zeta(y)=\lim \zeta(y^m)$$
weakly in  $L^q(\epsilon,T-\epsilon;L^q(\Gamma_1)) $ if $q<\infty$ and weakly in  $L^p(\epsilon,T-\epsilon;L^p(\Gamma_1)$ for all $2\le p<\infty$ if $q=\infty$. 
By the lower semicontinuity of the norm in $L^p$ we find from \eqref{contrad1} that $\zeta(y)=0$.

Now we use (H5) to see that both $y$ and $\displaystyle\frac{\partial y}{\partial n}$ are zero on $\Sigma_1$.  At this point we consider a slightly larger domain, $\tilde \Omega= \Omega\cup \{\Gamma_1+B_\delta(0)\}$. Consider $\tilde y$ the extension with $0$ of $y$ to $\tilde Q=(0,T)\times\tilde\Omega$ and $\tilde g$ also the extension with $0$ of $g$. Observe that $\tilde y$ is variational solution corresponding to $\tilde g$ and with null boundary conditions on the new piece of boundary. We are now in position to apply the strong maximum principle for weak solutions to parabolic systems obtained in \cite{lefter_melnig2024} we get that $\tilde y$ is zero in $\tilde Q$ and thus the source $\tilde g$ must be zero in $\tilde Q$, which  contradicts the fact that $g\not=0$.
 
\fin
\subsection*{Comparison with the results of O. Yu Imanuvilov and M. Yamamoto}

We mention here the result of O. Yu Imanuvilov and M. Yamamoto \cite{imayam1998} obtained in $L^2$ for the source of  one parabolic  equation. It is important to see that the shape of $\Omega\subset\R^N$ is not restricted and the homogeneous boundary conditions are prescribed on a part of the boundary $\Gamma_0$ and the observation is made on $\Gamma_1=\partial\Omega\setminus\Gamma_0$:

\begin{equation}\label{sysinitial2}\tag{I-Y}
\left\lbrace
\begin{array}{ll}
D_ty-\sum\limits_{j,k=1}^N D_j(a^{jk} D_ky)+\sum\limits_{k=1}^N b^{k} D_ky+cy=g &(0,T)\times\Omega,
\\
\beta(x)\frac{\partial y }{\partial n_{A}}+\eta(x)y=0  &(0,T)\times\Gamma_0, \\
\end{array} \quad
\right.\end{equation}

For  a fixed $\theta\in(0,T) $ an observation instant  of time, they have obtained  $L^2$ boundary  estimates in the  following class of sources:
\begin{equation}\label{Gset}
\mathcal{G}_{2,\tilde{c}}=\left\{
g\in W^{1,1}((0,T);L^2(\Omega)):\,\left|\frac{\partial g(t,x)}{\partial t}\right|\leq
\tilde c|g(\theta,x)|,  {a.e. } (t,x)\in Q\right\}.
\end{equation}
In this context, their result provides estimates in terms of  more measured quantities involving the gradient and the time derivative of the solution  on the observed boundary:
\begin{theorem}
Let $g\in \mathcal{G}_{2,\tilde{c}}$  and $y\in W^{2,1}_2(Q)$ a solution to \eqref{sysinitial2} corresponding to $g$.  Then
\begin{equation}\label{I-Y-est}
\|g\|_{L^2}\leq C(\|y(\theta,\cdot)\|_{W^2_2(\Omega)}+\|y\|_{L^2(\Sigma_1)}+\|\nabla D_ty\|_{L^2(\Sigma_1)}+ \|\nabla y\|_{L^2(\Sigma_1)}),
\end{equation}
where $\Sigma_1=(0,T)\times\Gamma_1$.
\end{theorem}

The result is  based on the two following Carleman estimates:
\begin{proposition}\label{lemaCarlemanclassic} For $g\in L^2(Q)$,
there exist constants $\lambda_0=\lambda_0(\Omega),$ $s_0=s_0(\Omega)$ such that, for any $\lambda\geq\lambda_0$, $ s\geq s_0$ and some $C=C(T,\Omega)$, the following inequality holds:
\begin{equation}\label{classicalCarleman}
\begin{aligned}
&\int_{Q}\left[(s\varphi)^{p-1}\left(|D_ty|^2+|D^2y|^2\right)+(s\varphi)^{p+1}|Dy|^2+(s\varphi)^{p+3}|y|^2\right]e^{2s\alpha}dxdt\\
&\leq C\int_{Q}(s\varphi)^p|g|^2e^{2s\alpha}dxdt+\\
&+C\int_{[0,T]\times\Gamma_1}\left((s\varphi)^{p}|D_t y|^2+(s\varphi)^{p+1}|\nabla y|^2
+(s\varphi)^{p+3}|y|^2\right)e^{2s\alpha}d\sigma,\quad p\in\{0,1\}
\end{aligned} 
\end{equation}
for $y\in H^1(0,T; L^2(\Omega))\cap L^2(0,T; H^2(\Omega))$ solution of \eqref{sysinitial2}.
\end{proposition}

Observe that if we impose homogeneous Dirichlet boundary conditions on \( \Gamma_1 \), the Carleman estimates by O. Yu. Imanuvilov and M. Yamamoto directly yield our Carleman estimate. This follows from the fact that, in this case, we have \( D_t y |_{\partial\Omega} = 0 \) and \( \nabla y |_{\partial\Omega} = \frac{\partial y}{\partial n} \) on \( \Gamma_1 \).  

As a result, our main conclusion remains valid for any domain, regardless of its topology, as long as homogeneous Dirichlet boundary conditions are imposed. However, this does not hold when Neumann or Robin boundary conditions are prescribed on \( \Gamma_1 \). In these cases, the estimate \eqref{I-Y-est} requires additional measurements involving the tangential and time derivatives of the solution on \( \Gamma_1 \), making the situation more restrictive.



\newpage


\begin{thebibliography}{10}

\bibitem{barbu2000}
V.~{Barbu}.
\newblock {Exact controllability of the superlinear heat equation.}
\newblock {\em {Appl. Math. Optim.}}, 42(1):73--89, 2000.

\bibitem{bar12002}
	Viorel {Barbu}.
	\newblock {The Carleman inequality for linear parabolic equations in $L^q$
		norm.}
	\newblock {\em {Differential Integral Equations}}, 15(5):513--525, 2002.
	

\bibitem{corgueros2010}
Jean-Michel {Coron}, Sergio {Guerrero}, and Lionel {Rosier}.
\newblock {Null controllability of a parabolic system with a cubic coupling
  term.}
\newblock {\em {SIAM J. Control Optim.}}, 48(8):5629--5653, 2010.

\bibitem{furima1996}
A.V. {Fursikov} and O.Yu. {Imanuvilov}.
\newblock {\em {Controllability of evolution equations.}}
\newblock Seoul: Seoul National Univ., 1996.

\bibitem{fernandez-Zuazua2000}
Enrique Fernández-Cara, Enrique Zuazua,
\newblock {\em {Null and approximate controllability for weakly blowing up semilinear heat equations.}}
\newblock {\em {Ann. Inst. H. Poincaré Anal. Non Linéaire}},   17 (2000), no. 5, pp. 583–616.

\bibitem{imayam1998}
Oleg~Yu {Imanuvilov} and Masahiro {Yamamoto}.
\newblock {Lipschitz stability in inverse parabolic problems by the Carleman
  estimate.}
\newblock {\em {Inverse Probl.}}, 14(5):1229--1245, 1998.

\bibitem{lady}
O.~A. Ladyzenskaja, V.~A. Solonnikov, and N.~N. Uralceva.
\newblock {\em Linear and quasilinear equations of parabolic type}.
\newblock Vol. 23. American Mathematical Society, Providence, R.I., 1968.

\bibitem{balch}
K\'{e}vin Le~Balc'h.
\newblock Controllability of a {$4\times 4$} quadratic reaction-diffusion
  system.
\newblock {\em J. Differential Equations}, 266(6):3100--3188, 2019.

\bibitem{lefter_melnig2024}
C.-G. Lefter and E.-A. Melnig.
\newblock Reaction-diffusion systems in annular domains: source stability estimates with boundary observations.
\newblock{\em {Optimization, 1–34}, https://doi.org/10.1080/02331934.2024.2444618.}


\bibitem{melnig2020}
E.-A. Melnig.
\newblock Stability in ${L}^q$-norm for inverse source parabolic problems.
\newblock {\em Journal of Inverse and Ill-posed Problems}, 28(6):797--814,
  2020.

\bibitem{melnig2021}
E.-A. Melnig.
\newblock {Stability in inverse source problems for nonlinear
  reaction-diffusion systems}.
\newblock {\em NoDEA Nonlinear Differential Equations Appl.}, 28(45), 2021.

\bibitem{wein2}
Murray~H. Protter and Hans~F. Weinberger.
\newblock {\em Maximum principles in differential equations}.
\newblock Springer-Verlag, New York, 1984.
\newblock Corrected reprint of the 1967 original.

\bibitem{wein}
Hans~F. Weinberger.
\newblock Invariant sets for weakly coupled parabolic and elliptic systems.
\newblock {\em Rend. Mat. (6)}, 8:295--310, 1975.
\newblock Collection of articles dedicated to Mauro Picone on the occasion of
  his ninetieth birthday.

\end{thebibliography}
\end{document}